\documentclass[pre,aps,twocolumn]{revtex4}
\usepackage{epsfig}
\usepackage{bm}

\begin{document}

\title{Correlations between large prime numbers}

\author{\small  A. Bershadskii}
\affiliation{\small {ICAR, P.O.B. 31155, Jerusalem 91000, Israel}}

\begin{abstract}
It is shown that short-range correlations between large prime numbers ($\sim 10^5$ and larger) 
have a Poissonian nature. Correlation length $\zeta \simeq 4.5$ for the primes $\sim 10^5$ and 
it is increasing logarithmically according to the prime number theorem. For moderate prime numbers 
($\sim 10^4$) the Poissonian distribution is not applicable (while the correlation length $\zeta$ 
surprisingly continues to follow to the logarithmical law. A chaotic (deterministic) hypothesis has 
been suggested to explain the moderate prime numbers apparent randomness.  \\


\end{abstract}

\maketitle

\begin{figure} \vspace{-0.5cm}\centering
\epsfig{width=.45\textwidth,file=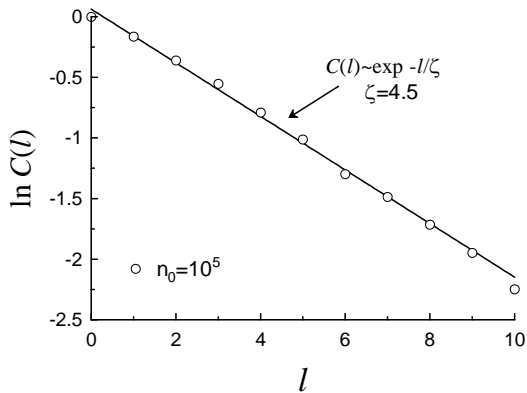} \vspace{-4.5cm}
\caption{Autocorrelation of the $v(n)$ function in a window with the width 
$\Delta =10^4$ centered at $n_0 \sim 10^5$.}
\end{figure}
\begin{figure} \vspace{-0.5cm}\centering
\epsfig{width=.45\textwidth,file=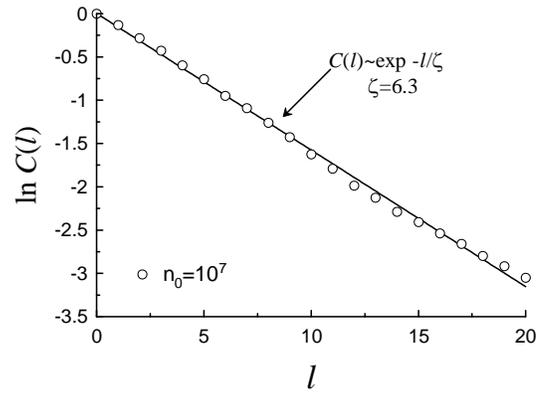} \vspace{-5cm}
\caption{Autocorrelation of the $v(n)$ function in a window with the width 
$\Delta =10^5$ centered at $n_0 \sim 10^7$.}
\end{figure}

{\it Introduction}.  The prime number distribution is apparently random. The apparent randomness, however, can be stochastic 
or chaotic (deterministic). The nature of the randomness is still unsolved problem. Using probabilistic 
methods we should take into account that there are different levels of averaging. One can consider the probabilities 
as representing a first level of averaging while the correlation functions represent a second one. Although 
the first level contains more details the second one can provide a more robust picture. In order to use 
the second level of averaging let us define a binary function $v(n)$ of integers $n=3,4,5,....$, which 
takes two values +1 or -1 and changes its sign passing any prime number. This function contains full 
information about the prime numbers distribution. Due to the prime number theorem, which states that the 
"average length" ($\tau_p$) of the gap between a prime $p$ and the next prime number is proportional 
(asymptotically) to $\ln p$, the $v(n)$ is a statistically non-stationary function. However, for sufficiently large $n$ this 
function can be considered (in a certain level of averaging) as a stationary function in the windows 
centered at $n_0$ with width $\Delta/n_0 \ll 1$. An heuristic estimate of a relation:
$$
\frac{\tau_{(p+\Delta)}}{\tau_p} \sim \frac{\ln (p +\Delta)}{\ln p} \simeq 1+\frac{\Delta}{p\ln p}  \eqno{(1)} 
$$
for $\Delta/p \ll 1$ provides a basis for such consideration. We will show (numerically) that 
on the level of the correlation function this idea indeed seems to be rather well applicable for sufficiently 
large prime numbers. On the level of the probabilities, however, application of this idea demands certain 
additional average (after that the results obtained by the two different methods can be compared). \\

{\it Large prime numbers}.   
Before performing this analysis let us recall certain simple properties of the random telegraph signal. 
The random telegraph signal is a binary process $x(n)$ which takes two values +1 or -1, and  has Poissonian 
distribution of the values of $n$ where $x(n)$ changes its sign. For the {\it stationary} case the Poisson distribution results in the exponential distribution of gaps $\tau$ between consecutive values of $n$ where $x(n)$ changes 
its sign
$$
P(\tau) = \frac{1}{\Theta} e^{-\tau/\Theta}  \eqno{(2)}
$$
with {\it constant} value of $\Theta \geq 0$. The autocorrelation function of the statistically stationary 
random telegraph signal 
$x(n)$
$$
C(l)= \langle x(n)x(n + l)\rangle - \langle x(n)\rangle \langle x(n + l)\rangle = e^{-l/\zeta} \eqno{(3)}
$$
where the angular brackets denote averaging over realizations and
$$
\zeta = \frac{\Theta}{2}   \eqno{(4)}
$$
\begin{figure} \vspace{-1cm}\centering
\epsfig{width=.45\textwidth,file=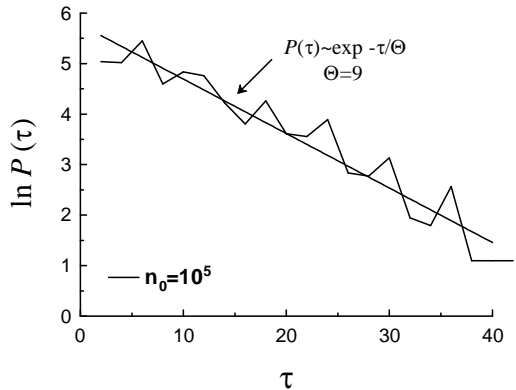} \vspace{-5cm}
\caption{Distribution of the gaps between consecutive primes for $v(n)$ in a 'window of stationarity' with the width $\Delta =10^4$ centered at $n_0 \sim 10^5$.}
\end{figure}
\begin{figure} \vspace{-0.2cm}\centering
\epsfig{width=.45\textwidth,file=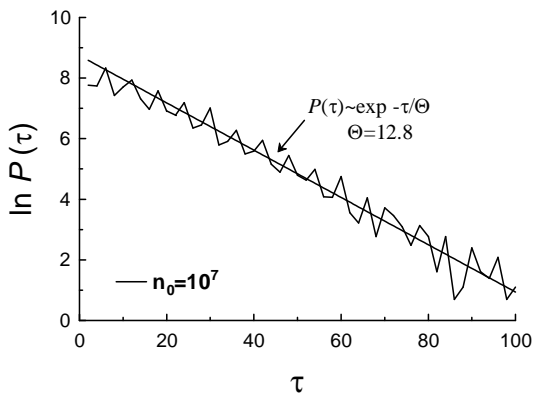} \vspace{-5cm}
\caption{As in Fig. 3 but in a window with the width 
$\Delta =10^5$ centered at $n_0 \sim 10^7$.}
\end{figure}
\begin{figure} \vspace{-0.5cm}\centering
\epsfig{width=.45\textwidth,file=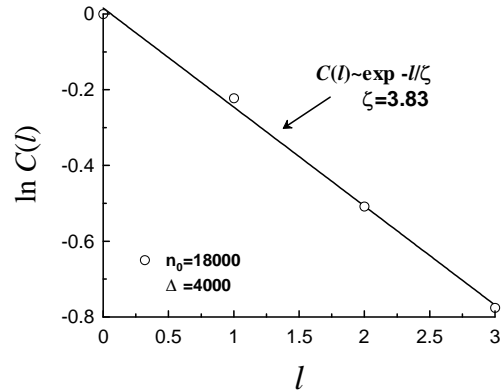} \vspace{-5.5cm}
\caption{Autocorrelation of the $v(n)$ function in a window with the width 
$\Delta =4000$ centered at $n_0 \sim 18000$.}
\end{figure}
 Figure 1 shows a correlation function computed for the $v(n)$ function in a 'window of stationarity' 
with the width $\Delta =10^4$ centered at $n_0 \sim 10^5$. The semi-log axes has been used in this figure 
in order to indicate exponential decay, Eq. (3), of the autocorrelation function (the straight line). 
The correlation length $\zeta \simeq 4.5$. Figure 2 shows analogous autocorrelation function but in  
a window with the width $\Delta =10^5$ centered at $n_0 \sim 10^7$. The correlation length $\zeta \simeq 6.3$ 
in this case. Figures 3 and 4 show the distribution of the gaps $\tau$ between consecutive primes: 
$P(\tau)$, for $v(n)$ in a 'window of stationarity' with the width $\Delta =10^4$ centered at $n_0 \sim 10^5$ 
and in a window with the width 
$\Delta =10^5$ centered at $n_0 \sim 10^7$ respectively. The straight lines in these figures represent the best 
fit to the data (also certain kind of averaging) and the slopes of these best fits provide us with the 
exponents of the Eq. (2): $\Theta \simeq  9$ and $\Theta \simeq  12.8$ respectively. Now one can compare 
the values of the correlation length $\zeta$ with values of the $\Theta$ and with the Poissonian relation 
Eq. (4). One can see that in the second level of averaging correspondence to the Poissonian distribution of the 
primes in the considered windows is rather good. Similarly to the "average gap" (see above) the second level's parameters $\zeta$ and $\Theta=2\zeta$ are proportional to $\ln n_0$ due to the prime number theorem. \\

{\it Moderate prime numbers}. Figure 5 shows an exponentially decaying part of the autocorrelation computed 
for the $v(n)$ function in a 'window of stationarity' with the width $\Delta =4000$ centered at $n_0 \sim 18000$ 
(moderate values of the primes). Naturally, the exponential range here is rather short. However, the 
computed value of the correlation length $\zeta \simeq 3.83$ follows quite precisely to the logarithmic dependence on $n_0$. It is surprising if we take into account that the distribution of the gaps between consecutive 
primes (the parameter $\Theta$) does not, as one can see from figure 6 (the $\Theta \simeq 9$ here). 
This means, in particular, that the $v(n)$ function in this 
window is not a random telegraph signal. Something wrong is going with the Poissonian hypothesis for these values of the 
primes. But if previously we related the logarithmical dependence of the correlation length $\zeta$ on $n_0$ 
just with the logarithmical dependence of $\Theta$, then how can we explain it now? 

The problem with the Poissonian distribution for the moderate primes may be related to the 
strong and persistent periodic oscillations in the distribution of the gaps. Figure 7 shows spectrum of these 
oscillations (seen in Fig. 6) after a detrending, and the strong peak in Fig. 7 corresponds to the period 
equal to 6 (as well as for the larger $n_0$). In the Refs. \cite{wol}-\cite{ac} this period has been studied 
both numerically and theoretically. We are more interested, however, in nature of the randomness of the 
moderate primes. The spectrum of the function $v(n)$ in the moderate window of stationarity, shown in 
figure 8 in the semi-log scales, may help in solving this problem. One can see a pronounced range of an 
exponential decay of the spectrum. Many of the well known chaotic attractors ('Lorenz', 'R\"{o}ssler', etc.) exhibit the exponentially decaying 
spectra and this spectrum is considered as a strong indication of a chaotic behavior \cite{sig}-\cite{fm}. 
\begin{figure} \vspace{-1cm}\centering
\epsfig{width=.45\textwidth,file=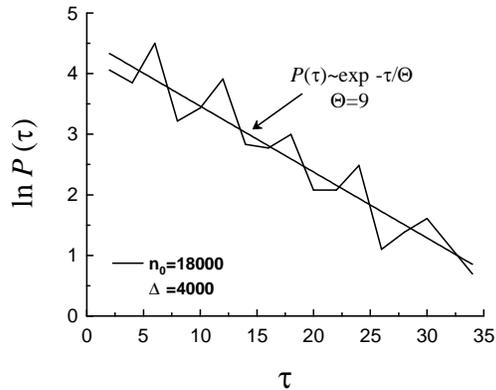} \vspace{-5cm}
\caption{Distribution of the gaps between consecutive primes for $v(n)$ in a 'window of stationarity' with the width $\Delta =4000$ centered at $n_0 \sim 18000$.}
\end{figure}
\begin{figure} \vspace{-0.3cm}\centering
\epsfig{width=.45\textwidth,file=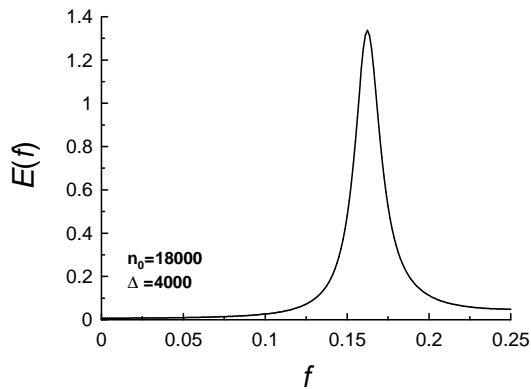} \vspace{-5cm}
\caption{Spectrum of the oscillations shown in Fig. 6 (after a detrending). }
\end{figure}
\begin{figure} \vspace{-0.5cm}\centering
\epsfig{width=.45\textwidth,file=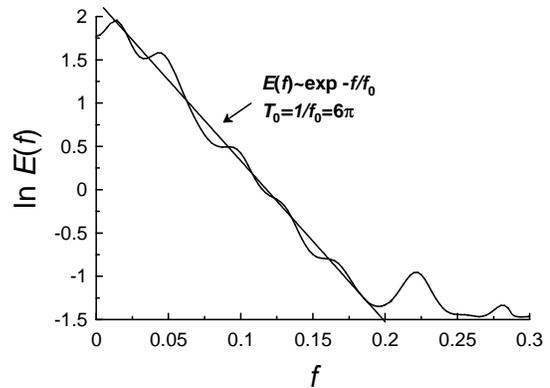} \vspace{-5.5cm}
\caption{Spectrum of the function $v(n)$ in a 'window of stationarity' with the 
width $\Delta =4000$ centered at $n_0 \sim 18000$.}
\end{figure}
Here, for comparison with the Fig. 8, we will consider a chaotic spectrum generated by 
the Kaplan-Yorke map \cite{ky}. In the Langevin 
approach to Brownian motion the equation of motion is 
$$
\dot{y}= -\gamma y + F(t)  \eqno{(5)}
$$
where the fluctuating kick force on the particle is a Gaussian white noise: 
$F(t)= \sum_n \eta_n \delta (t-n\tau)$ and $y(t)$ to take values in $R^m$.
One can assume \cite{beck} that the evolution of the kick strengths 
is determined by a discrete time dynamical system $T$ on the phase space and projected onto $R^m$ 
by a function $g$:
$$
\eta_n = g(x_{n-1}),~~~~  x_{n+1}= T x_n  \eqno{(6)}
$$
Then the solution of Eq. (5) is
$$
y(t) = e^{-\gamma(t-n\tau)} y_n \eqno{(7)}
$$
where $n$ equals the integer part of the relation $t/\tau$ and the recurrence
$$
x_{n+1}= T x_n, ~~~ y_{n+1}=\alpha y_n + g(x_n)  \eqno{(8)}
$$
provides value of $y_n$ (with $\alpha = e^{-\gamma\tau}$). In certain sense the dynamical
system (8) is equivalent to the stochastic differential equation (5). In the 
generalization related to the Eq. (8) the force $F(t)$ can be considered as a {\it non}-Gaussian 
process which is determined by $g$ and $T$. 

The Kaplan-Yorke map \cite{ky},\cite{jo},\cite{spot} is a particular 
simple case for this generalization:
$$
Tx = 2x~ (mod~ 1), ~~~ g(x)= \cos 4 \pi x  \eqno{(9)}
$$
An optimal computational algorithm for the Kaplan-Yorke map is 
$$
a_{n+1}=2a_n~mod(p), ~~x_{n+1}=a_n/p, ~y_{n+1}=\alpha y_n + \cos 4 \pi x \eqno{(10)}
$$
where $p$ is a large prime number.

Figure 9 shows spectrum of a chaotic solution of the Kaplan-Yorke map ($\alpha =0.2$). We used the semi-log 
axes in this figure in order to indicate exponential decay of the spectrum (the straight line). 

  One can compare Fig. 9 and Fig. 8 in order to conclude on the chaotic nature of the "window defined" distribution of the moderate primes. 
\begin{figure} \vspace{-1cm}\centering
\epsfig{width=.45\textwidth,file=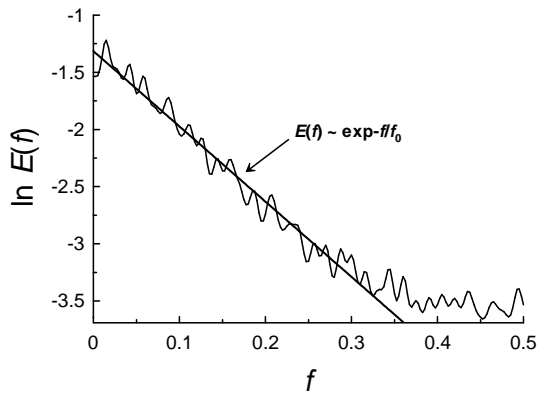} \vspace{-5cm}
\caption{Spectrum of a chaotic solution of the Kaplan-Yorke map. }
\end{figure}
 
   A rigorously minded reader can find a lot of good mathematics related the above discussed approach: such as Gallagher theorem for short intervals \cite{gal} (about Poisson distribution in the limit case of logarithmically short windows), 
Cramer's probabilistic model \cite{cra} and the difficulties of this model in the short intervals (Maier's theorem 
\cite{ma},\cite{gra},\cite{ms}). Since present consideration was devoted to the correlations between large primes 
and the correlation length $\zeta$ is rather small (see Figs. 1,2 and 5) the consideration was restricted by the 
windows of statistic stationarity where the Poisson distribution seems to be reigning on the second level of statistical averaging at sufficiently large primes. Above these windows distributions of the primes can be different (see, for instance, Ref. \cite{sou} and references therein). For moderate values of primes the Poisson distribution 
seems to be not applicable and a chaotic hypothesis has been suggested instead.


\begin{thebibliography}{99}
\bibitem{wol} M. Wolf, Applications of statistical mechanics in
number theory, Physica A, {\bf 274}, 149-157 (1999)
\bibitem{kum} P. Kumar, P.Ch. Ivanov, and H.E. Stanley, Information Entropy and Correlations in Prime Numbers, 
arXiv:cond-mat/0303110v4 [cond-mat.stat-mech], (2003).
\bibitem{ac} S. Aresa, and M. Castro, Hidden structure in the randomness
of the prime number sequence?, Physica A, {\bf 360}, 285-296 (2006).
\bibitem{sig} D.E. Sigeti, Survival of deterministic dynamics in the presence
of noise and the exponential decay of power spectra at high frequency, 
Phys. Rev. E, {\bf 52}, 2443; Physica D, {\bf 82}, 136 (1995).
\bibitem{o} N. Ohtomo, K. Tokiwano, Y. Tanaka et. al., Exponential Characteristics
of Power Spectral Densities Caused by Chaotic Phenomena, J. Phys. Soc. Jpn., {\bf 64} 1104 (1995). 
\bibitem{fa} J. D. Farmer, Chaotic Attractors of an Infinite-Dimensional Dynamical
System, Physica D, {\bf 4}, 366 (1982).
\bibitem{fm} U. Frisch and R. Morf, Intermittency in nonlinear dynamics and
singularities at complex times, Phys. Rev., {\bf 23}, 2673 (1981).
\bibitem{ky} J.L. Kaplan and J.A. Yorke, in Functional Differential Equations and Approximations of Fixed Points, 
Lecture Notes in Mathematics, {\bf 730}, p. 204. (New York: Springer 1979).
\bibitem{beck} C. Beck, Ergodic Properties of a Kicked Damped Particle, Commun. Math. Phys., {\bf 130}, 51 (1990).
\bibitem{jo} R.V. Jensen and C.R. Oberman, Calculation of the Statistical Properties of Strange Attractors,
Phys. Rev. Lett., {\bf 46}, 1547 (1981).
\bibitem{spot} J. C. Sprott, Chaos and Time-Series Analysis (Oxford Univ. Press, 2003)
\bibitem{gal} P.X. Gallagher, On the distribution of primes in short intervals, Mathematika, {\bf 23}, 
4-9, (1976).
\bibitem{cra} H.Cramer, On the order of magnitude of the difference between consecutive prime numbers. Acta
Arith., {\bf 2}, 23-46 (1936).
\bibitem{ma} H. Maier, Primes in short intervals, The Michigan Math. Journal, {\bf 32} 221 (1985).
\bibitem{gra} A. Granville, Proc. Intern. Congress of Math., {\bf I},  388 (1995) 
(Zurich, Switzerland, 1994).
\bibitem{ms}H.L. Montgomery, and K. Soundararajan, Primes in short intervals, Commun. Math. Phys., {\bf 252}, 589-617 (2004).
\bibitem{sou} K. Soundararajan, The distribution of the prime numbers, 
NATO Science Series II: Mathematics, Physics and Chemistry, , {\bf 237}, 59-83 (2007).
\end{thebibliography}
\end{document}